\documentclass[10pt, A4paper]{article}

\usepackage{graphicx,url}
\usepackage{amssymb,amsmath}
\usepackage{epstopdf}
\usepackage{amssymb,amsmath,amsthm}

\usepackage{setspace}

\usepackage{url,color}
\usepackage{lscape}

\numberwithin{equation}{section}
\date{25 January 2010}

\begin{document}

\title{Equations of Mathematical Physics and Compositions of Brownian and
Cauchy processes}
\author{L.Beghin\thanks{%
\textquotedblleft Sapienza\textquotedblright\ University of Rome}, L.Sakhno%
\thanks{%
Kyiv National Taras Shevchenko University}, E.Orsingher\thanks{%
\textquotedblleft Sapienza\textquotedblright\ University of Rome.
Corresponding author}}
\maketitle

\begin{abstract}
We consider different types of processes obtained by composing Brownian
motion $B(t)$, fractional Brownian motion $B_{H}(t)$ and Cauchy processes $%
C(t)$ in different manners.

We study also multidimensional iterated processes in $\mathbb{R}^{d},$ like,
for example, $\left( B_{1}(|C(t)|),...,B_{d}(|C(t)|)\right) $ and $\left(
C_{1}(|C(t)|),...,C_{d}(|C(t)|)\right) ,$ deriving the corresponding partial
differential equations satisfied by their joint distribution.

We show that many important partial differential equations, like wave
equation, equation of vibration of rods, higher-order heat equation, are
satisfied by the laws of the iterated processes considered in the work.

Similarly we prove that some processes like $%
C(|B_{1}(|B_{2}(...|B_{n+1}(t)|...)|)|)$ are governed by fractional
diffusion equations.

\


\end{abstract}

\section{Introduction}

The iterated Brownian motion is obtained by the composition of two
independent Brownian motions $B_{1}$ and $B_{2},$ as follows:%
\begin{equation}
I(t)=B_{1}(|B_{2}(t)|),\qquad t>0.  \label{due.1}
\end{equation}%
Recently this kind of processes has been studied by many authors: see, for
example, Burdzy (1994), Khoshnevisan and Lewis (1996), Allouba and Zheng
(2001), Allouba (2002), Beghin and Orsingher (2009). As far as the
applications are concerned, it has been observed that the iterated Brownian
motion $I$ is suitable to describe diffusions in cracks (De Blassie (2004))
and many other physical phenomena. In particular, fractures on a rectangular
slab can be viewed at as trajectories of a Brownian motion (see Chudnovsky
and Kunin (1987)). The flow of a gas on a fracture can be represented by a
Brownian motion moving on a Brownian sample path and therefore this model
produces an iterated Brownian motion (see, for example, Khoshnevisan and
Lewis (1999)).

In Orsingher and Zhao (1999) it is stated that the density of the process%
\begin{equation}
p_{I}(x,t)=2\int_{0}^{\infty }\frac{e^{-\frac{x^{2}}{2s}}}{\sqrt{2\pi s}}%
\frac{e^{-\frac{s^{2}}{2t}}}{\sqrt{2\pi t}}ds  \label{due.2}
\end{equation}%
satisfies the fourth-order equation%
\begin{equation}
\frac{\partial p}{\partial t}(x,t)=\frac{1}{2^{3}}\frac{\partial ^{4}p}{%
\partial x^{4}}+\frac{1}{2\sqrt{2\pi t}}\frac{d^{2}}{dx^{2}}\delta (x),
\label{due.4}
\end{equation}%
for $x\in \mathbb{R}\;$and $t>0.$

In De Blassie (2004) a similar result is proved under the initial condition $%
p(x,0)=f(x)$.

Moreover it has been shown in Orsingher and Beghin (2004) that the density (%
\ref{due.2}) is a solution to the fractional equation%
\begin{equation}
\frac{\partial ^{\frac{1}{2}}p}{\partial t^{\frac{1}{2}}}=\frac{1}{2^{3/2}}%
\frac{\partial ^{2}p}{\partial x^{2}},\qquad x\in \mathbb{R},\;t>0
\label{due.3}
\end{equation}%
(with initial condition $p(x,0)=\delta (x)$), which is a limiting case of
the fractional telegraph equation%
\begin{equation*}
\frac{\partial p}{\partial t}+2^{3/2}\lambda \frac{\partial ^{\frac{1}{2}}p}{%
\partial t^{\frac{1}{2}}}=c^{2}\frac{\partial ^{2}p}{\partial x^{2}},\qquad
x\in \mathbb{R},\;c,\lambda ,t>0,
\end{equation*}%
for $\lambda ,c\rightarrow \infty ,$ in such a way that $c^{2}/\lambda
\rightarrow 1.$

Generalizations of this result have been obtained in many directions. First
of all by considering the $n$-times iterated Brownian motion%
\begin{equation}
I_{n}(t)=B_{1}(|B_{2}(...|B_{n+1}(t)|...)|),  \label{due.5}
\end{equation}%
involving $n+1$ independent one-dimensional Brownian motions.

It has been shown in Orsingher and Beghin (2009) that the law of (\ref{due.5}%
) satisfies the following fractional equation%
\begin{equation}
\frac{\partial ^{\frac{1}{2^{n}}}p}{\partial t^{\frac{1}{2^{n}}}}=2^{\frac{1%
}{2^{n}}-2}\frac{\partial ^{2}p}{\partial x^{2}},\qquad x\in \mathbb{R}%
,\;t>0,  \label{due.6}
\end{equation}%
with $p(x,0)=\delta (x).$

For the vector process
\begin{equation}
I^{d}(t)=\left\{
\begin{array}{c}
B_{1}(|B(t)|) \\
... \\
... \\
B_{d}(|B(t)|)%
\end{array}%
\right. ,\qquad t>0,  \label{id}
\end{equation}%
where $B$ and $B_{1},...,B_{d}$ are mutually independent Brownian motions,
it is proved in Orsingher and Beghin (2009) that the joint law $%
p_{I}^{d}=p_{I}^{d}(x_{1},x_{2},...,x_{d},t)$ is a solution to the
fractional equation%
\begin{equation}
\frac{\partial ^{\frac{1}{2}}p}{\partial t^{\frac{1}{2}}}=\frac{1}{2^{3/2}}%
\sum_{j=1}^{d}\frac{\partial ^{2}p}{\partial x_{j}^{2}},\qquad x_{j}\in
\mathbb{R},\mathbb{\;}j=1,...,d,\;t>0.  \label{due.9}
\end{equation}%
The constants appearing in equations (\ref{due.3}), (\ref{due.6}) and (\ref%
{due.9}) have been chosen in such a way that the Brownian motions involved
in the construction of the iterated processes possess standard distribution.

If we consider the slightly more general fractional equation%
\begin{equation}
\frac{\partial ^{\frac{1}{2}}p}{\partial t^{\frac{1}{2}}}=\lambda
^{2}\sum_{j=1}^{d}\frac{\partial ^{2}p}{\partial x_{j}^{2}},\qquad x_{j}\in
\mathbb{R},\mathbb{\;}j=1,...,d,\;t>0,  \label{due.10}
\end{equation}%
then the joint density of the vector process reads%
\begin{equation}
p_{I}^{d}(x_{1},x_{2},...,x_{d},t)=2\int_{0}^{\infty }\frac{e^{-\frac{w^{2}}{%
2(2^{3}\lambda ^{4}t)}}}{\sqrt{2\pi (2^{3}\lambda ^{4}t)}}\prod_{k=1}^{d}%
\frac{e^{-\frac{x_{k}^{2}}{2w}}}{\sqrt{2\pi w}}dw.  \label{due.11}
\end{equation}

Result (\ref{due.11}) shows that the components of the $d$-dimensional
vector process $\left( B_{1}(|B(t)|),...,B_{d}(|B(t)|)\right) $\ are no
longer independent and the parameter $\lambda $ enters into the variance of
the "time process" $B(t),t>0.$

For the law of the vector process $I^{d}\left( t\right) $\ we can also write
the equation which is analogous to (\ref{due.4}):
\begin{equation}
\frac{\partial p}{\partial t}(x,t)=\frac{1}{2^{3}}\left( \sum_{j=1}^{d}\frac{%
\partial ^{2}}{\partial x_{j}^{2}}\right) ^{2}p+\frac{1}{2\sqrt{2\pi t}}%
\sum_{j=1}^{d}\frac{d^{2}}{dx_{j}^{2}}\prod_{k=1}^{d}\delta (x_{k}),
\label{due.44}
\end{equation}%
with initial condition $p(x_{1},...,x_{d},0)=\prod_{k=1}^{d}\delta (x_{k}).$

As far as the $d$-dimensional case is concerned, along with the definition (%
\ref{id}) of the process $I^{d}$\ with components being independent
one-dimensional iterated Brownian motions, another choice is to define the
process as $B_{\left( d\right) }(|B(t)|),t>0$,\ with $B_{\left( d\right) }$\
being a $d$-dimensional Brownian motion independent of the Brownian motion $%
B(t).$\ For the process $I^{d}(t)$\ Allouba and Zheng (2001) and De Blassie
(2004) show that the function $%
u(t,x)=E_{x}(f(I^{d}(t)))=E(f(I^{d}(t))|I_{0}^{d}=x)$\ solves the initial
value problem
\begin{equation}
\frac{\partial }{\partial t}u(t,x)=\Delta ^{2}u(t,x)+\frac{1}{\sqrt{\pi t}}%
\Delta f(x),\text{ \ \ }u(0,x)=f(x),  \label{deb}
\end{equation}%
for $t>0$, $x\in \mathbb{R}^{d}$, and `sufficiently good' functions $f(x)$, $%
x\in \mathbb{R}^{d}$\ (which are supposed to be bounded with bounded Hölder
continuous second derivatives).

Moreover, when replacing the outer process $B_{\left( d\right) }$\ with a
continuous Markov process $X$\ and considering $Z(t)=X(|B(t)|)$, equation (%
\ref{deb}) remains valid for $u(t,x)=E_{x}(f(Z(t)))$\ if we replace the
Laplacian with the generator $L_{x}$\ of the continuous semigroup associated
with the Markov process $X$, that is:
\begin{equation*}
\frac{\partial }{\partial t}u(t,x)=L_{x}^{2}u(t,x)+\frac{1}{\sqrt{\pi t}}%
L_{x}f(x),\text{ \ \ }u(0,x)=f(x),
\end{equation*}%
$f\in D(L_{x})$, where $D(L_{x})$ is the domain of the operator $L_{x}.$\
Moreover $u(t,x)=E_{x}(f(Z(t)))$\ solves the fractional Cauchy problem
\begin{equation*}
\frac{\partial ^{\frac{1}{2}}}{\partial t^{\frac{1}{2}}}u(t,x)=L_{x}u(t,x),%
\text{ \ }u(0,x)=f(x),\text{ \ }x\in \mathbb{R}^{d},\;t>0
\end{equation*}%
(see, Baeumer, Meerschaert and Nane (2009)).

A further extension is given in Baeumer, Meerschaert and Nane (2009) (see
also Nane (2008)). Let $X\left( t\right) =x+X_{0}\left( t\right) ,$\ where $%
X_{0}\left( t\right) $\ is a Lévy process in $\mathbb{R}^{d}$\ starting at
zero, and let $L_{x}$\ be the generator of the semigroup $%
p(x,t)=E_{x}(f(X(t)))$. Then for any $f\in D(L_{x})$\ and any $n=2,3,4,...$\
it was shown that both Cauchy problems for the higher-order partial
differential equation
\begin{equation}
\frac{\partial u}{\partial t}=\sum_{i=1}^{n-1}\frac{t^{\frac{j}{n}-1}}{%
\Gamma \left( \frac{j}{n}\right) }L_{x}^{j}f+L_{x}^{n}u,\text{ \ }%
u(x,0)=f(x),\qquad x\in \mathbb{R}^{d},t>0  \label{due.7}
\end{equation}%
and for the fractional equation
\begin{equation}
\frac{\partial ^{\frac{1}{n}}u}{\partial t^{\frac{1}{n}}}=L_{x}u,\text{ \ \ }%
u(x,0)=f(x),\text{ \ \ }x\in \mathbb{R}^{d},\;t>0,  \label{due.8}
\end{equation}%
have the same unique solution. This is given by $%
u(x,t)=E_{x}(f(X(E^{1/n}(t))))$\ $,$\ where $E^{1/n}(t)$\ is the inverse Lé%
vy subordinator of order $\frac{1}{n},$\ that is, $E^{1/n}(t)=\inf
\{s:D(s)>t\},$\ where $D(s)$\ is a stable subordinator with index $\frac{1}{n%
}$.

For the Markov process with the inverse Lévy subordinator of order $\frac{1}{%
2}$\ one gets the same governing equations as in the case of a Markov
process with a Brownian time, that is, these two processes have the same
one-dimensional distributions.

In the light of the above equivalence of (\ref{due.7}) and (\ref{due.8})
(and inspired by the result of Orsingher and Beghin (2009), equation (\ref%
{due.6})), the following PDE connection of $n$-iterated Brownian motion
appears in Nane (2008): the unique solution of the Cauchy problems (\ref%
{due.7}) and (\ref{due.8}) with $n=2^{k}$\ is given by $%
u(x,t)=E_{x}(f(X(I_{k-1}(t)))),$\ where $X$\ is the Lévy process described
above and $I_{k}(t)$\ is $k$-times iterated Brownian motion (\ref{due.5}).

For the law of the $n$-times iterated Brownian motion, we can write an
analog of equations (\ref{due.4}) and (\ref{due.44}):
\begin{equation}
\frac{\partial p}{\partial t}=\sum_{i=1}^{m-1}\frac{t^{\frac{j}{m}-1}}{%
\Gamma \left( \frac{j}{m}\right) }\left( \frac{\partial ^{2}}{\partial x^{2}}%
\right) ^{j}\delta (x)+\left( \frac{\partial ^{2}}{\partial x^{2}}\right)
^{m}p,\text{ \ }p(x,0)=\delta (x),\text{ \ }x\in \mathbb{R},\text{ \ }t>0,
\label{due.444}
\end{equation}%
where $m=2^{n}.$

For the fractional diffusion equation of order $1/n,$ $n\in \mathbb{N},$%
\begin{equation}
\frac{\partial ^{\frac{1}{n}}p}{\partial t^{\frac{1}{n}}}=\frac{1}{2}\frac{%
\partial ^{2}p}{\partial x^{2}},\qquad x\in \mathbb{R},\;t>0,  \label{due.12}
\end{equation}%
with initial condition $p_{1/n}(x,t)=\delta (x),$ it has been shown in
Beghin and Orsingher (2003) that the solution coincides with the density of
the process%
\begin{equation}
I_{1/n}(t)=B\left( \prod_{j=1}^{n-1}G_{j}(t)\right) ,\qquad n>1,
\label{due.13}
\end{equation}%
where $\left( G_{1}(t),...,G_{n}(t)\right) $ is a random variable with joint
distribution%
\begin{equation}
p(w_{1},...,w_{n-1})=\frac{n^{\frac{n-1}{2}}}{(2\pi )^{\frac{n-1}{2}}\sqrt{t}%
}e^{-\frac{w_{1}^{n}+...+w_{n-1}^{n}}{\sqrt[n-1]{n^{n}t}}%
}w_{2}...w_{n-1}^{n-2},\qquad w_{j}\geq 0,\;1\leq j\leq n-1.  \label{due.14}
\end{equation}

Taking into account the result of Baeumer, Meerschaert and Nane (2009) on
the solution to the equation (\ref{due.8}), we can conclude that the process
(\ref{due.13}) and the process $B\left( E^{1/n}(t)\right) ,$ that is a
Brownian motion subordinated to the inverse Lévy process of order $\frac{1}{n%
}$, have the same one-dimensional distributions, which can be written down
in a closed form.

In D'Ovidio and Orsingher (2009) various types of processes obtained by
composing independent fractional Brownian motions $B_{H_{1}}(t)$ and $%
B_{H_{2}}(t),t>0$ (with Hurst parameters respectively equal to $%
H_{1},H_{2}\in \left( 0,1\right) $) have been studied. In particular, for
the density $p_{H_{1}H_{2}}(x,t)$ of $B_{H_{1}}(|B_{H_{2}}(t)|)$ it has been
shown that the governing equation has the following structure:%
\begin{equation*}
(1+H_{1}H_{2})t\frac{\partial p}{\partial t}+t^{2}\frac{\partial ^{2}p}{%
\partial t^{2}}=H_{1}^{2}H_{2}^{2}\left\{ 2x\frac{\partial p}{\partial x}%
+x^{2}\frac{\partial ^{2}p}{\partial x^{2}}\right\} ,\quad x\in \mathbb{R}%
,t>0.
\end{equation*}%
In this paper we extend some of these results in several directions, by
considering different compositions of Brownian motions, fractional Brownian
motions and Cauchy processes.

We start by studying the iterated Brownian motion defined in (\ref{due.1})
in the case where $B_{1}$ is endowed by drift, giving the two governing
equations, the fractional and fourth-order one (see (\ref{gri}) below).

Another direction of our research concerns processes obtained by combining
Brownian motions (possibly fracional Brownian motions) and Cauchy processes.
In particular we study the multidimensional vector processes $\left(
B_{1}(|C(t)|),...,B_{d}(|C(t)|)\right) $ and $\left(
C_{1}(|B(t)|),...,C_{d}(|B(t)|)\right) ,$ which involve independent Brownian
motions $B_{j}(t),t>0$ and Cauchy processes $C_{j}(t),t>0.$

For $d=1$ it is has been proved in D'Ovidio and Orsingher (2009) that the
law of $B(|C(t)|),t>0$, given by%
\begin{equation}
p_{BC}(x,t)=\frac{2}{\pi }\int_{0}^{\infty }\frac{e^{-\frac{x^{2}}{2s}}}{%
\sqrt{2\pi s}}\frac{t}{t^{2}+s^{2}}ds,\quad x\in \mathbb{R},\text{ }t>0
\end{equation}%
satisfies the following fourth-order equation%
\begin{equation}
\frac{\partial ^{2}p}{\partial t^{2}}=-\frac{1}{2^{2}}\frac{\partial ^{4}p}{%
\partial x^{4}}-\frac{1}{\pi t}\frac{d^{2}}{dx^{2}}\delta (x),\qquad x\in
\mathbb{R},\text{ }t>0,
\end{equation}%
with initial condition $p(x,0)=\delta (x).$ We extend this result to the
case where the Brownian motion is substituted by the fractional Brownian
motion $B_{H}(t),t>0,$ with Hurst parameter $H\in \left( 0,1\right) .$ We
show that the density function of $B_{H}(|C(t)|),t>0$ resolves the following
equation%
\begin{equation*}
t^{2}\frac{\partial ^{2}p}{\partial t^{2}}=-\left[ H(H-1)\frac{\partial }{%
\partial x}x-H^{2}\frac{\partial ^{2}}{\partial x^{2}}x^{2}\right] p-\frac{%
2Ht}{\pi }\frac{d^{2}}{dx^{2}}\delta (x)1_{\left\{ 0<H\leq \frac{1}{2}%
\right\} },
\end{equation*}%
with $p(x,0)=\delta (x).$

Also the case of a Cauchy process with a random time represented by an
iterated Brownian motion is taken into account and the corresponding
fractional governing equation is derived.

The iterated Cauchy process, defined as $J_{CC}(t)=C_{1}(|C_{2}(t)|),t>0$
has been considered in D'Ovidio and Orsingher (2009) and its connection with
the wave equation
\begin{equation}
\frac{\partial ^{2}p}{\partial t^{2}}=\frac{\partial ^{2}p}{\partial x^{2}}-%
\frac{1}{\pi tx^{2}},\qquad x\in \mathbb{R},\;t>0
\end{equation}%
has been established.

We will show here that when the process $C_{1}$ is endowed with a position
parameter $a\neq 0$, the distribution of the iterated Cauchy process
satisfies the following wave equation%
\begin{equation}
\frac{\partial ^{2}p}{\partial t^{2}}=\frac{\partial ^{2}p}{\partial x^{2}}-%
\frac{1}{\pi t(x-a)^{2}},\qquad x\in \mathbb{R},\;t>0.
\end{equation}%
Finally we consider the $d$-dimensional vector $\left(
C_{1}(|C(t)|),...,C_{d}(|C(t)|)\right) $ and obtain the equation governed by
its density, i.e.%
\begin{equation}
p_{CC}^{d}(x_{1},...,x_{d},t)=\frac{2}{\pi ^{d+1}}\int_{0}^{+\infty
}\prod\limits_{j=1}^{d}\frac{s}{s^{2}+x_{j}^{2}}\frac{t}{t^{2}+s^{2}}ds,
\end{equation}%
where the use of partial Riesz fractional derivatives is required.

\section{Iterated Brownian motion with drift}

If the processes composing the iterated Brownian motions possess drift, the
fractional equations and the higher-order equations governing the
distributions are somewhat different. We start by considering the process $%
I^{\mu }(t)=B_{1}^{\mu }(|B_{2}(t)|),t>0$, with law%
\begin{equation}
p_{I}^{\mu }(x,t)=2\int_{0}^{\infty }\frac{e^{-\frac{(x-\mu s)^{2}}{2s}}}{%
\sqrt{2\pi s}}\frac{e^{-\frac{s^{2}}{2t}}}{\sqrt{2\pi t}}ds.  \label{due.15}
\end{equation}

It has been shown in Beghin and Orsingher (2009) that (\ref{due.15}) solves
the following fractional equation of order $\nu =1/2:$%
\begin{equation}
\frac{\partial ^{\frac{1}{2}}p}{\partial t^{\frac{1}{2}}}=\frac{1}{2^{3/2}}%
\frac{\partial ^{2}p}{\partial x^{2}}-\frac{\mu }{\sqrt{2}}\frac{\partial p}{%
\partial x},\qquad x\in \mathbb{R},\;t>0.  \label{due.16}
\end{equation}%
As a check we evaluate the Laplace-Fourier transform of (\ref{due.15}):%
\begin{eqnarray}
&&\int_{0}^{\infty }e^{-\eta t}\int_{-\infty }^{+\infty }e^{i\beta x}\left(
2\int_{0}^{\infty }\frac{e^{-\frac{(x-\mu s)^{2}}{2s}}}{\sqrt{2\pi s}}\frac{%
e^{-\frac{s^{2}}{2t}}}{\sqrt{2\pi t}}ds\right) dxdt  \label{due.17} \\
&=&2\int_{0}^{\infty }e^{-\frac{1}{2}\beta ^{2}s+i\beta \mu s}\frac{e^{-s%
\sqrt{2\eta }}}{\sqrt{2\eta }}ds=\frac{\eta ^{\frac{1}{2}-1}}{\frac{\beta
^{2}}{2^{3/2}}-\frac{i\beta \mu }{\sqrt{2}}+\sqrt{\eta }}.  \notag
\end{eqnarray}

The previous result coincides with the Fourier-Laplace transform of (\ref%
{due.16}).

We show in the following theorem that, as in the case where $\mu =0,$ the
density of $I^{\mu }(t)$ satisfies also a fourth-order equation.

\

\noindent \textbf{Theorem 2.1 }\textit{The density of the process }$I^{\mu
}(t)=B_{1}^{\mu }(|B_{2}(t)|),t>0$\textit{\ given in (\ref{due.15}) is a
solution to the following equation}%
\begin{equation}
\frac{\partial }{\partial t}p(x,t)=\frac{1}{2}\left( \frac{1}{2}\frac{%
\partial ^{2}}{\partial x^{2}}-\mu \frac{\partial }{\partial x}\right)
^{2}p(x,t)+\frac{1}{\sqrt{2\pi t}}\left( \frac{1}{2}\frac{d^{2}}{dx^{2}}-\mu
\frac{d}{dx}\right) \delta (x),\qquad x\in \mathbb{R},\;t>0  \label{gri}
\end{equation}%
\textit{with initial condition }$p(x,0)=\delta (x).$

\noindent \textbf{Proof \ }By taking the first-order derivative of (\ref%
{due.15}), we have that%
\begin{eqnarray}
&&\frac{\partial }{\partial t}p_{I}^{\mu }(x,t)  \label{due.20} \\
&=&2\int_{0}^{\infty }\frac{e^{-\frac{(x-\mu s)^{2}}{2s}}}{\sqrt{2\pi s}}%
\frac{\partial }{\partial t}\frac{e^{-\frac{s^{2}}{2t}}}{\sqrt{2\pi t}}ds
\notag \\
&=&\int_{0}^{\infty }\frac{e^{-\frac{(x-\mu s)^{2}}{2s}}}{\sqrt{2\pi s}}%
\frac{\partial ^{2}}{\partial s^{2}}\frac{e^{-\frac{s^{2}}{2t}}}{\sqrt{2\pi t%
}}ds  \notag \\
&=&\left[ \text{integrating by parts}\right]  \notag \\
&=&-\int_{0}^{\infty }\frac{\partial }{\partial s}\frac{e^{-\frac{(x-\mu
s)^{2}}{2s}}}{\sqrt{2\pi s}}\frac{\partial }{\partial s}\frac{e^{-\frac{s^{2}%
}{2t}}}{\sqrt{2\pi t}}ds  \notag \\
&=&-\left. \frac{\partial }{\partial s}\frac{e^{-\frac{(x-\mu s)^{2}}{2s}}}{%
\sqrt{2\pi s}}\frac{e^{-\frac{s^{2}}{2t}}}{\sqrt{2\pi t}}\right\vert
_{0}^{\infty }+\int_{0}^{\infty }\frac{\partial ^{2}}{\partial s^{2}}\frac{%
e^{-\frac{(x-\mu s)^{2}}{2s}}}{\sqrt{2\pi s}}\frac{e^{-\frac{s^{2}}{2t}}}{%
\sqrt{2\pi t}}ds  \notag \\
&=&\frac{1}{\sqrt{2\pi t}}\left( \frac{1}{2}\frac{d^{2}}{dx^{2}}-\mu \frac{d%
}{dx}\right) \delta (x)+\int_{0}^{\infty }\frac{\partial }{\partial s}\left(
\frac{1}{2}\frac{\partial ^{2}}{\partial x^{2}}-\mu \frac{\partial }{%
\partial x}\right) \frac{e^{-\frac{(x-\mu s)^{2}}{2s}}}{\sqrt{2\pi s}}\frac{%
e^{-\frac{s^{2}}{2t}}}{\sqrt{2\pi t}}ds  \notag \\
&=&\frac{1}{2}\left( \frac{1}{2}\frac{\partial ^{2}}{\partial x^{2}}-\mu
\frac{\partial }{\partial x}\right) ^{2}p_{I}^{\mu }(x,t)+\frac{1}{\sqrt{%
2\pi t}}\left( \frac{1}{2}\frac{d^{2}}{dx^{2}}-\mu \frac{d}{dx}\right)
\delta (x).  \notag
\end{eqnarray}%
\hfill $\square $

\

\noindent \textbf{Remark 2.1 }The differential operator appearing in the
fourth-order equation (\ref{gri}) is the formal square of the operator
appearing in the fractional equation (\ref{due.16}). In the special case
where $\mu =0$, we obtain again equation (\ref{due.4}).

\

We consider now the case where the process representing the
\textquotedblleft time" possesses drift, then the law of the iterated
process $I^{\mu }(t)=B_{1}(|B_{2}^{\mu }(t)|),t>0$, reads%
\begin{equation}
q_{I}^{\mu }(x,t)=\frac{1}{C(t)}\int_{0}^{\infty }\frac{e^{-\frac{x^{2}}{2s}}%
}{\sqrt{2\pi s}}\frac{e^{-\frac{(s-\mu t)^{2}}{2t}}}{\sqrt{2\pi t}}ds,
\end{equation}%
where
\begin{equation*}
C(t)=\Pr \left\{ B_{2}^{\mu }(t)>0\right\} =\int_{-\mu \sqrt{t}}^{\infty }%
\frac{e^{-\frac{y^{2}}{2}}}{\sqrt{2\pi }}dy
\end{equation*}%
is the normalizing factor for the density of $|B_{2}^{\mu }(t)|.$

The Fourier transform of $q_{I}^{\mu }(x,t)$ becomes
\begin{equation}
\int_{-\infty }^{+\infty }e^{i\beta x}q_{I}^{\mu }(x,t)dx=\frac{1}{C(t)}%
\int_{-\mu \sqrt{t}}^{\infty }\frac{e^{-\frac{\beta ^{2}}{2}(\sqrt{t}w+\mu
t)-\frac{w^{2}}{2}}}{\sqrt{2\pi }}dw.  \label{due.19}
\end{equation}%
The evaluation of the Laplace transform of (\ref{due.19}) poses serious
problems. For this reason the case where the process representing the
\textquotedblleft time" possesses drift is not further developed here.

\section{Iterated processes involving the Cauchy process}

The transition function of a centered Cauchy process $C(t),t>0$ is given by%
\begin{equation}
p_{C}(x,t)=\frac{t}{\pi (t^{2}+x^{2})},\quad x\in \mathbb{R},\text{ }t>0
\label{qua.1}
\end{equation}%
and it is well-known that $p_{C}$ is a solution to the Laplace equation%
\begin{equation}
\left( \frac{\partial ^{2}}{\partial x^{2}}+\frac{\partial ^{2}}{\partial
t^{2}}\right) p_{C}(x,t)=0.  \label{qua.2}
\end{equation}

Furthermore $p_{C}$ is a solution to the following space-fractional equation%
\begin{equation}
\frac{\partial }{\partial t}p_{C}(x,t)=-\frac{\partial }{\partial |x|}%
p_{C}(x,t),  \label{qua.3}
\end{equation}%
where
\begin{equation}
\frac{\partial }{\partial |x|}f(x)=\frac{1}{\pi }\frac{d}{dx}\left[
\int_{-\infty }^{x}\frac{f(y)}{x-y}dy-\int_{x}^{\infty }\frac{f(y)}{y-x}dy%
\right]  \label{qua.4}
\end{equation}%
is a special case (for $\nu =1$) of the Riesz fractional space-derivative
and possesses Fourier transform equal to%
\begin{equation}
\int_{-\infty }^{+\infty }e^{i\beta x}\frac{\partial }{\partial |x|}%
f(x)dx=|\beta |\int_{-\infty }^{+\infty }e^{i\beta x}f(x)dx=|\beta |\mathcal{%
F(}\beta \mathcal{)}.  \label{fou}
\end{equation}

By composing the standard Brownian motion and the Cauchy process, we obtain
the following new processes:%
\begin{equation}
J_{BC}(t)=B(|C(t)|),\qquad t>0  \label{qua.5}
\end{equation}%
and
\begin{equation}
J_{CB}(t)=C(|B(t)|),\qquad t>0,  \label{qua.6}
\end{equation}%
where $C$ and $B$ are mutually independent.

An alternative definition of the first-order fractional derivative is given
in Saichev and Zaslavsky (1997) (see formula (A.39)) and reads%
\begin{equation}
\frac{d}{d|x|}f(x)=-\frac{1}{\pi }\int_{0}^{+\infty }\frac{%
f(x-y)-2f(x)+f(x+y)}{y^{2}}dy.  \label{qua.7}
\end{equation}%
This is a special case of%
\begin{equation}
\frac{d^{\nu }}{d|x|^{\nu }}f(x)=\frac{1}{2\Gamma (-\nu )\cos \frac{\pi \nu
}{2}}\int_{0}^{+\infty }\frac{f(x-y)-2f(x)+f(x+y)}{y^{2}}dy,  \label{qua.8}
\end{equation}%
for $0<\nu \leq 1$. We note that (\ref{qua.8}) reduces to (\ref{qua.7}) for $%
\nu =1,$ because, by the reflection formula of the Gamma function, the
constant can be rewritten as
\begin{eqnarray*}
\frac{1}{2\Gamma (-\nu )\cos \frac{\pi \nu }{2}} &=&-\frac{1}{2\cos \frac{%
\pi \nu }{2}}\frac{\Gamma (1+\nu )}{\frac{\pi }{\sin \pi \nu }} \\
&=&-\frac{\Gamma (1+\nu )\sin \frac{\pi \nu }{2}}{\pi },
\end{eqnarray*}%
which, for $\nu =1$, yields $-1/\pi .$

We can convert definition (\ref{qua.4}) into (\ref{qua.7}) by means of
integrations by parts, after suitable changes of variables.

As a further check we can show that the Fourier transform of (\ref{qua.7})
coincides with (\ref{fou}):%
\begin{eqnarray}
&&-\frac{1}{\pi }\int_{-\infty }^{+\infty }e^{i\beta x}\left(
\int_{0}^{+\infty }\frac{f(x-y)-2f(x)+f(x+y)}{y^{2}}dy\right) dx
\label{quaqua} \\
&=&-\frac{1}{\pi }\int_{-\infty }^{+\infty }f(w)e^{i\beta
w}dw\int_{0}^{+\infty }\frac{e^{i\beta y}}{y^{2}}dy+\frac{2}{\pi }%
\int_{-\infty }^{+\infty }f(x)e^{i\beta x}dx\int_{0}^{+\infty }\frac{1}{y^{2}%
}dy+  \notag \\
&&-\frac{1}{\pi }\int_{-\infty }^{+\infty }f(w)e^{i\beta
w}dw\int_{0}^{+\infty }\frac{e^{-i\beta y}}{y^{2}}dy  \notag \\
&=&-\frac{1}{\pi }\mathcal{F}(\beta )\int_{0}^{+\infty }\frac{e^{i\beta
y}-2+e^{-i\beta y}}{y^{2}}dy.  \notag
\end{eqnarray}%
The last integral can be evaluated as follows:%
\begin{eqnarray*}
&&\int_{0}^{+\infty }\frac{e^{i\beta y}-1}{y^{2}}dy+\int_{0}^{+\infty }\frac{%
e^{-i\beta y}-1}{y^{2}}dy \\
&=&\int_{0}^{+\infty }\frac{1}{y^{2}}dy\int_{0}^{y}i\beta e^{i\beta
w}dw-\int_{0}^{+\infty }\frac{1}{y^{2}}dy\int_{0}^{y}i\beta e^{-i\beta w}dw
\\
&=&i\beta \int_{0}^{+\infty }e^{i\beta w}dw\int_{w}^{+\infty }\frac{1}{y^{2}}%
dy-i\beta \int_{0}^{+\infty }e^{-i\beta w}dw\int_{w}^{+\infty }\frac{1}{y^{2}%
}dy \\
&=&-\left[ 2\beta \int_{0}^{+\infty }\frac{e^{i\beta w}}{2wi}dw-2\beta
\int_{0}^{+\infty }\frac{e^{-i\beta w}}{2wi}dw\right] \\
&=&-2\beta \int_{0}^{+\infty }\frac{\sin \beta w}{w}dw=-\pi |\beta |,
\end{eqnarray*}%
which, inserted into (\ref{quaqua}), gives (\ref{fou}).

It has been shown in D'Ovidio and Orsingher (2009) that the transition
density of the process $J_{BC}(t)=B(|C(t)|),$ $t>0$, which is given by%
\begin{equation}
p_{BC}(x,t)=\frac{2}{\pi }\int_{0}^{\infty }\frac{e^{-\frac{x^{2}}{2s}}}{%
\sqrt{2\pi s}}\frac{t}{t^{2}+s^{2}}ds,\quad x\in \mathbb{R},\text{ }t>0
\label{qua.9}
\end{equation}%
satisfies the following equation%
\begin{equation}
\frac{\partial ^{2}p}{\partial t^{2}}=-\frac{1}{2^{2}}\frac{\partial ^{4}p}{%
\partial x^{4}}-\frac{1}{\pi t}\frac{d^{2}}{dx^{2}}\delta (x),\qquad x\in
\mathbb{R},\text{ }t>0,  \label{qua.10}
\end{equation}%
with initial condition $p(x,0)=\delta (x).$

\

\noindent \textbf{Remark 3.1 }In mathematical physics the equation%
\begin{equation}
\frac{\partial ^{2}u}{\partial t^{2}}=-K\frac{\partial ^{4}u}{\partial x^{4}}
\label{qua.9bis}
\end{equation}%
represents the vibration of rods (see Elmore and Heald (1969), p.116).
Equation (\ref{qua.9bis}) coincides with (\ref{qua.10}) for $x\neq 0.$

\

\noindent \textbf{Remark 3.2 }The density (\ref{qua.9}) can be worked out in
an equivalent form, by applying the subordinating relationship%
\begin{equation*}
\frac{t}{\pi (t^{2}+s^{2})}=\int_{0}^{\infty }\frac{e^{-\frac{s^{2}}{2w}}}{%
\sqrt{2\pi w}}\frac{te^{-\frac{t^{2}}{2w}}}{\sqrt{2\pi w^{3}}}dw.
\end{equation*}%
Thus the transition density of the process $I_{BC}(t)=B(|C(t)|),$ $t>0$ can
be rewritten as%
\begin{equation}
p_{BC}(x,t)=2\int_{0}^{\infty }\frac{e^{-\frac{x^{2}}{2s}}}{\sqrt{2\pi s}}%
\left( \int_{0}^{\infty }\frac{e^{-\frac{s^{2}}{2w}}}{\sqrt{2\pi w}}\frac{%
te^{-\frac{t^{2}}{2w}}}{\sqrt{2\pi w^{3}}}dw\right) ds.  \label{new}
\end{equation}%
Formula (\ref{new}) corresponds to the law of an iterated Brownian motion $%
I(t)=B_{1}(|B_{2}(t)|),$ taken at a random time $T(t),$ which coincides in
distribution with the first-passage time of a standard Brownian motion. Thus
(\ref{new}) coincides with the density of the process $I(T(t)),t>0,$ where $%
T(t)=\inf (s:B(s)=t).$

\

We consider now the composition of a fractional Brownian motion with a
Cauchy process, $J_{B_{H}C}(t)=B_{H}(|C(t)|),t>0$.

\

\noindent \textbf{Theorem 3.1 }\textit{The density of }$J_{B_{H}C}(t),t>0$%
\textit{, which is given by}%
\begin{equation}
p_{B_{H}C}(x,t)=\frac{2}{\pi }\int_{0}^{+\infty }\frac{e^{-\frac{x^{2}}{%
2s^{2H}}}}{\sqrt{2\pi s^{2H}}}\frac{t}{t^{2}+s^{2}}ds,  \label{bc}
\end{equation}%
\textit{satisfies the following equation}%
\begin{equation}
t^{2}\frac{\partial ^{2}p}{\partial t^{2}}=-\left[ H(H-1)\frac{\partial }{%
\partial x}x-H^{2}\frac{\partial ^{2}}{\partial x^{2}}x^{2}\right] p-\frac{%
2Ht}{\pi }\frac{d^{2}}{dx^{2}}\delta (x)1_{\left\{ 0<H\leq \frac{1}{2}%
\right\} },  \label{bc2}
\end{equation}%
\textit{with initial condition }$p(x,0)=\delta (x).$

\noindent \textbf{Proof }We take the second order time-derivative of (\ref%
{bc}), that is%
\begin{eqnarray}
\frac{\partial ^{2}}{\partial t^{2}}p_{B_{H}C}(x,t) &=&\frac{2}{\pi }%
\int_{0}^{+\infty }\frac{e^{-\frac{x^{2}}{2s^{2H}}}}{\sqrt{2\pi s^{2H}}}%
\frac{\partial ^{2}}{\partial t^{2}}\left( \frac{t}{t^{2}+s^{2}}\right) ds
\label{bc3} \\
&=&-\frac{2}{\pi }\int_{0}^{+\infty }\frac{e^{-\frac{x^{2}}{2s^{2H}}}}{\sqrt{%
2\pi s^{2H}}}\frac{\partial ^{2}}{\partial s^{2}}\left( \frac{t}{t^{2}+s^{2}}%
\right) ds  \notag \\
&=&\left. \frac{2}{\pi }\frac{\partial }{\partial s}\left( \frac{e^{-\frac{%
x^{2}}{2s^{2H}}}}{\sqrt{2\pi s^{2H}}}\right) \frac{t}{t^{2}+s^{2}}%
\right\vert _{s=0}^{s=+\infty }+  \notag \\
&&-\frac{2}{\pi }\int_{0}^{+\infty }\frac{\partial ^{2}}{\partial s^{2}}%
\left( \frac{e^{-\frac{x^{2}}{2s^{2H}}}}{\sqrt{2\pi s^{2H}}}\right) \frac{t}{%
t^{2}+s^{2}}ds.  \notag
\end{eqnarray}%
The first term in (\ref{bc3}) can be worked out by considering that the
transition density of the fractional Brownian motion $p_{H}=p_{H}(x,t)$
satisfies the following equation%
\begin{equation*}
\frac{\partial p}{\partial t}=Ht^{2H-1}\frac{\partial ^{2}p}{\partial x^{2}}
\end{equation*}%
and it becomes%
\begin{equation}
-\left. \frac{2H}{\pi }\frac{ts^{2H-1}}{t^{2}+s^{2}}\frac{\partial ^{2}}{%
\partial x^{2}}\left( \frac{e^{-\frac{x^{2}}{2s^{2H}}}}{\sqrt{2\pi s^{2H}}}%
\right) \right\vert _{s=0}=\left\{
\begin{array}{l}
0\qquad \frac{1}{2}<H<1 \\
-\frac{2Ht}{\pi }\frac{d^{2}}{dx^{2}}\delta (x)\qquad 0<H\leq \frac{1}{2}%
\end{array}%
\right. .  \label{bc9}
\end{equation}%
The second term in (\ref{bc3}) can be evaluated as follows%
\begin{eqnarray}
&&-\frac{2}{\pi }\int_{0}^{+\infty }\frac{\partial }{\partial s}\left(
Hs^{2H-1}\frac{\partial ^{2}}{\partial x^{2}}\left( \frac{e^{-\frac{x^{2}}{%
2s^{2H}}}}{\sqrt{2\pi s^{2H}}}\right) \right) \frac{t}{t^{2}+s^{2}}ds
\label{bc4} \\
&=&-\frac{2H(2H-1)}{\pi }\int_{0}^{+\infty }s^{2H-2}\frac{\partial ^{2}}{%
\partial x^{2}}\left( \frac{e^{-\frac{x^{2}}{2s^{2H}}}}{\sqrt{2\pi s^{2H}}}%
\right) \frac{t}{t^{2}+s^{2}}ds+  \notag \\
&&-\frac{2H^{2}}{\pi }\int_{0}^{+\infty }s^{4H-2}\frac{\partial ^{4}}{%
\partial x^{4}}\left( \frac{e^{-\frac{x^{2}}{2s^{2H}}}}{\sqrt{2\pi s^{2H}}}%
\right) \frac{t}{t^{2}+s^{2}}ds  \notag \\
&=&+\frac{2H(2H-1)}{\pi }\frac{\partial }{\partial x}\left(
x\int_{0}^{+\infty }s^{-2}\frac{e^{-\frac{x^{2}}{2s^{2H}}}}{\sqrt{2\pi s^{2H}%
}}\frac{t}{t^{2}+s^{2}}ds\right) +  \notag \\
&&+\frac{2H^{2}}{\pi }\int_{0}^{+\infty }s^{2H-2}\frac{\partial ^{3}}{%
\partial x^{3}}\left( x\frac{e^{-\frac{x^{2}}{2s^{2H}}}}{\sqrt{2\pi s^{2H}}}%
\right) \frac{t}{t^{2}+s^{2}}ds.  \notag
\end{eqnarray}%
We rewrite the last line in (\ref{bc4}) as%
\begin{eqnarray}
&&\frac{2H^{2}}{\pi }\int_{0}^{+\infty }s^{2H-2}\frac{\partial ^{2}}{%
\partial x^{2}}\left( \frac{e^{-\frac{x^{2}}{2s^{2H}}}}{\sqrt{2\pi s^{2H}}}%
\right) \frac{t}{t^{2}+s^{2}}ds+  \label{bc5} \\
&&-\frac{2H^{2}}{\pi }\int_{0}^{+\infty }s^{-2}\frac{\partial ^{2}}{\partial
x^{2}}\left( x^{2}\frac{e^{-\frac{x^{2}}{2s^{2H}}}}{\sqrt{2\pi s^{2H}}}%
\right) \frac{t}{t^{2}+s^{2}}ds  \notag \\
&=&-\frac{2H^{2}}{\pi }\int_{0}^{+\infty }s^{-2}\frac{\partial }{\partial x}%
\left( x\frac{e^{-\frac{x^{2}}{2s^{2H}}}}{\sqrt{2\pi s^{2H}}}\right) \frac{t%
}{t^{2}+s^{2}}ds+  \notag \\
&&-\frac{2H^{2}}{\pi }\int_{0}^{+\infty }s^{-2}\frac{\partial ^{2}}{\partial
x^{2}}\left( x^{2}\frac{e^{-\frac{x^{2}}{2s^{2H}}}}{\sqrt{2\pi s^{2H}}}%
\right) \frac{t}{t^{2}+s^{2}}ds.  \notag
\end{eqnarray}%
By inserting (\ref{bc5}) into (\ref{bc4}) the second term of (\ref{bc3})
takes the form%
\begin{equation}
\left[ \frac{2H(H-1)}{\pi }\frac{\partial }{\partial x}x-\frac{2H^{2}}{\pi }%
\frac{\partial ^{2}}{\partial x^{2}}x^{2}\right] \int_{0}^{+\infty }s^{-2}%
\frac{e^{-\frac{x^{2}}{2s^{2H}}}}{\sqrt{2\pi s^{2H}}}\frac{t}{t^{2}+s^{2}}ds.
\label{bc6}
\end{equation}%
The integral appearing in (\ref{bc6}) can be calculated by observing that%
\begin{eqnarray}
&&\frac{1}{t}\int_{0}^{+\infty }s^{-2}\frac{e^{-\frac{x^{2}}{2s^{2H}}}}{%
\sqrt{2\pi s^{2H}}}\frac{t^{2}+s^{2}-s^{2}}{t^{2}+s^{2}}ds  \label{bc7} \\
&=&-\frac{1}{t^{2}}\frac{\pi }{2}p_{B_{H}C}(x,t)+\frac{1}{t}%
\int_{0}^{+\infty }s^{-2}\frac{e^{-\frac{x^{2}}{2s^{2H}}}}{\sqrt{2\pi s^{2H}}%
}ds.  \notag
\end{eqnarray}%
The integral in (\ref{bc7}) can be explicitly calculated as%
\begin{eqnarray*}
&&\int_{0}^{+\infty }s^{-2}\frac{e^{-\frac{x^{2}}{2s^{2H}}}}{\sqrt{2\pi
s^{2H}}}ds \\
&=&\left[ s=\frac{x^{1/H}}{2^{1/2H}}y^{-1/2H}\right] \\
&=&\int_{0}^{+\infty }\frac{e^{-y}x^{-2/H}y^{1/H}}{2^{-1/H}\sqrt{2\pi }}%
\frac{x^{1/H}}{2^{1/2H}}\left( -\frac{1}{2H}\right) y^{-\frac{1}{2H}-1}\frac{%
(2y)^{1/2}}{x}dy \\
&=&\frac{2^{\frac{1}{2H}-1}}{x^{\frac{1}{H}+1}H\sqrt{\pi }}\Gamma \left(
\frac{1}{2H}+\frac{1}{2}\right) ,
\end{eqnarray*}%
which, inserted in (\ref{bc7}), permits us to rewrite (\ref{bc6}) as follows:%
\begin{eqnarray}
&&\left[ \frac{2H(H-1)}{\pi }\frac{\partial }{\partial x}x-\frac{2H^{2}}{\pi
}\frac{\partial ^{2}}{\partial x^{2}}x^{2}\right] \left( -\frac{1}{t^{2}}%
\frac{\pi }{2}p_{B_{H}C}+\frac{2^{\frac{1}{2H}-1}}{x^{\frac{1}{H}+1}H\sqrt{%
\pi }t}\Gamma \left( \frac{1}{2H}+\frac{1}{2}\right) \right)  \notag \\
&=&-\frac{1}{t^{2}}\left[ H(H-1)\frac{\partial }{\partial x}x-H^{2}\frac{%
\partial ^{2}}{\partial x^{2}}x^{2}\right] p_{B_{H}C}+  \label{bc8} \\
&&+\frac{2^{\frac{1}{2H}}(H-1)}{\pi \sqrt{\pi }t}\Gamma \left( \frac{1}{2H}+%
\frac{1}{2}\right) \frac{d}{dx}\left( x^{-\frac{1}{H}}\right) -\frac{2^{%
\frac{1}{2H}}H}{\pi \sqrt{\pi }t}\Gamma \left( \frac{1}{2H}+\frac{1}{2}%
\right) \frac{d^{2}}{dx^{2}}\left( x^{-\frac{1}{H}+1}\right)  \notag \\
&=&-\frac{1}{t^{2}}\left[ H(H-1)\frac{\partial }{\partial x}x-H^{2}\frac{%
\partial ^{2}}{\partial x^{2}}x^{2}\right] p_{B_{H}C}+  \notag \\
&&-\frac{2^{\frac{1}{2H}}(H-1)}{\pi \sqrt{\pi }tH}\Gamma \left( \frac{1}{2H}+%
\frac{1}{2}\right) x^{-\frac{1}{H}-1}+\frac{2^{\frac{1}{2H}}(H-1)}{\pi \sqrt{%
\pi }tH}\Gamma \left( \frac{1}{2H}+\frac{1}{2}\right) x^{-\frac{1}{H}-1}
\notag \\
&=&-\frac{1}{t^{2}}\left[ H(H-1)\frac{\partial }{\partial x}x-H^{2}\frac{%
\partial ^{2}}{\partial x^{2}}x^{2}\right] p_{B_{H}C}.  \notag
\end{eqnarray}%
By putting together (\ref{bc8}) and (\ref{bc9}) we obtain equation (\ref{bc2}%
).\hfill$\square $

\

\noindent \textbf{Remark 3.3 }In order to check that for $H=\frac{1}{2}$
equation (\ref{bc2}) can be converted into (\ref{qua.10}), it is necessary
to perform the following calculations:%
\begin{eqnarray}
&&\frac{\partial ^{4}}{\partial x^{4}}\left( \frac{2}{\pi }\int_{0}^{+\infty
}\frac{e^{-\frac{x^{2}}{2s}}}{\sqrt{2\pi s}}\frac{t}{t^{2}+s^{2}}ds\right)
\label{bc10} \\
&=&-\frac{\partial ^{3}}{\partial x^{3}}\left( \frac{2}{\pi }%
x\int_{0}^{+\infty }\frac{e^{-\frac{x^{2}}{2s}}}{s\sqrt{2\pi s}}\frac{t}{%
t^{2}+s^{2}}ds\right)  \notag \\
&=&-\frac{2}{\pi }\frac{\partial ^{2}}{\partial x^{2}}\left(
\int_{0}^{+\infty }\frac{e^{-\frac{x^{2}}{2s}}}{s\sqrt{2\pi s}}\frac{t}{%
t^{2}+s^{2}}ds-x^{2}\int_{0}^{+\infty }\frac{1}{s^{2}}\frac{e^{-\frac{x^{2}}{%
2s}}}{\sqrt{2\pi s}}\frac{t}{t^{2}+s^{2}}ds\right) .  \notag
\end{eqnarray}%
In analogy with (\ref{bc7}), we can write%
\begin{equation}
\int_{0}^{+\infty }\frac{1}{s^{2}}\frac{e^{-\frac{x^{2}}{2s}}}{\sqrt{2\pi s}}%
\frac{t}{t^{2}+s^{2}}ds=-\frac{1}{t^{2}}\frac{\pi }{2}p_{BC}(x,t)+\frac{1}{t}%
\int_{0}^{+\infty }\frac{1}{s^{2}}\frac{e^{-\frac{x^{2}}{2s}}}{\sqrt{2\pi s}}%
ds.  \label{bc11}
\end{equation}%
Furthermore%
\begin{eqnarray}
&&\frac{\partial ^{2}}{\partial x^{2}}\left( \int_{0}^{+\infty }\frac{e^{-%
\frac{x^{2}}{2s}}}{s\sqrt{2\pi s}}\frac{t}{t^{2}+s^{2}}ds\right)
\label{bc12} \\
&=&-\frac{\partial }{\partial x}\left( x\int_{0}^{+\infty }\frac{e^{-\frac{%
x^{2}}{2s}}}{s^{2}\sqrt{2\pi s}}\frac{t}{t^{2}+s^{2}}ds\right)  \notag \\
&=&-\frac{\partial }{\partial x}\left[ x\left( -\frac{1}{t^{2}}\frac{\pi }{2}%
p_{BC}(x,t)+\frac{1}{t}\int_{0}^{+\infty }\frac{1}{s^{2}}\frac{e^{-\frac{%
x^{2}}{2s}}}{\sqrt{2\pi s}}ds\right) \right] .  \notag
\end{eqnarray}%
By inserting (\ref{bc11}) and (\ref{bc12}) into (\ref{bc10}) we have that%
\begin{eqnarray*}
\frac{\partial ^{4}}{\partial x^{4}}p_{BC}(x,t) &=&-\frac{1}{t^{2}}\frac{%
\partial ^{2}}{\partial x^{2}}\left( x^{2}p_{BC}(x,t)\right) -\frac{1}{t^{2}}%
\frac{\partial }{\partial x}\left( x\,p_{BC}(x,t)\right) + \\
&&+\frac{2}{\pi t}\frac{\partial ^{2}}{\partial x^{2}}\left(
x^{2}\int_{0}^{+\infty }\frac{e^{-\frac{x^{2}}{2s}}}{s^{2}\sqrt{2\pi s}}%
ds\right) +\frac{2}{\pi t}\frac{\partial }{\partial x}\left(
x\int_{0}^{+\infty }\frac{e^{-\frac{x^{2}}{2s}}}{s^{2}\sqrt{2\pi s}}ds\right)
\\
&=&-\frac{1}{t^{2}}\frac{\partial ^{2}}{\partial x^{2}}\left(
x^{2}p_{BC}(x,t)\right) -\frac{1}{t^{2}}\frac{\partial }{\partial x}\left(
x\,p_{BC}(x,t)\right) + \\
&&+\frac{2}{\pi t}\left[ \frac{\partial ^{2}}{\partial x^{2}}x^{2}+\frac{%
\partial }{\partial x}x\right] \int_{0}^{+\infty }\frac{e^{-\frac{x^{2}}{2s}}%
}{s^{2}\sqrt{2\pi s}}ds.
\end{eqnarray*}%
The calculations carried out in Theorem 3.1 permit us to conclude that, for $%
x\neq 0$,%
\begin{equation*}
\frac{\partial ^{4}}{\partial x^{4}}p_{BC}(x,t)=-\frac{1}{t^{2}}\frac{%
\partial ^{2}}{\partial x^{2}}\left( x^{2}p_{BC}(x,t)\right) -\frac{1}{t^{2}}%
\frac{\partial }{\partial x}\left( x\,p_{BC}(x,t)\right) =-\frac{\partial
^{2}}{\partial t^{2}}p_{BC}(x,t).
\end{equation*}%
Therefore, despite the formal diversity of (\ref{qua.10}) and (\ref{bc2}),
the law of $B(|C(t)|)$ satisfies also equation (\ref{bc2}).

\

We consider now the $d$-dimensional generalization of the results given in (%
\ref{qua.9}) and (\ref{qua.10}): let us define the vector process%
\begin{equation*}
J_{BC}^{d}(t)=\left\{
\begin{array}{c}
B_{1}(|C(t)|) \\
... \\
... \\
B_{d}(|C(t)|)%
\end{array}%
\right. ,\qquad t>0,
\end{equation*}%
where $B_{1},...,B_{d}$ are Brownian motions independent from each other and
from the Cauchy process $C(t),t>0.$ Then the following result holds.

\

\noindent \textbf{Theorem 3.2}\textit{\ The joint density of }$%
J_{BC}^{d}(t),t>0$\textit{, which is given by}%
\begin{equation}
p_{BC}^{d}(x_{1},...,x_{d},t)=\frac{2}{\pi }\int_{0}^{+\infty
}\prod_{j=1}^{d}\frac{e^{-\frac{x_{j}^{2}}{2s}}}{\sqrt{2\pi s}}\frac{t}{%
t^{2}+s^{2}}ds,  \label{qua.12}
\end{equation}%
\textit{is a solution to the following equation}%
\begin{equation}
\frac{\partial ^{2}p}{\partial t^{2}}=-\frac{1}{2^{2}}\left( \sum_{j=1}^{d}%
\frac{\partial ^{2}}{\partial x_{j}^{2}}\right) ^{2}p(x_{1},...,x_{d},t)-%
\frac{1}{\pi t}\sum_{j=1}^{d}\delta (x_{1})...\frac{d^{2}}{dx_{j}^{2}}\delta
(x_{j})...\delta (x_{d}),  \label{qua.13}
\end{equation}%
\textit{with initial condition }$p(x_{1},...,x_{d},0)=\prod_{j=1}^{d}\delta
(x_{j}).$

\noindent \textbf{Proof }By taking the second order time-derivative of (\ref%
{qua.12}) we get, by means of two integrations by parts, that%
\begin{eqnarray}
&&\frac{\partial ^{2}}{\partial t^{2}}p_{BC}^{d}(x_{1},...,x_{d},t)
\label{qua.14} \\
&=&-\frac{2}{\pi }\int_{0}^{+\infty }\frac{\partial ^{2}}{\partial s^{2}}%
\left( \prod_{j=1}^{d}\frac{e^{-\frac{x_{j}^{2}}{2s}}}{\sqrt{2\pi s}}%
\right) \frac{t}{t^{2}+s^{2}}ds+  \notag \\
&&+\left. \frac{2}{\pi }\frac{\partial }{\partial s}\left( \prod_{j=1}^{d}%
\frac{e^{-\frac{x_{j}^{2}}{2s}}}{\sqrt{2\pi s}}\right) \frac{t}{t^{2}+s^{2}}%
\right\vert _{0}^{+\infty }  \notag \\
&=&-\frac{2}{\pi }\left( \frac{1}{2}\sum_{j=1}^{d}\frac{\partial ^{2}}{%
\partial x_{j}^{2}}\right) ^{2}\int_{0}^{+\infty }\prod_{j=1}^{d}\frac{e^{-%
\frac{x_{j}^{2}}{2s}}}{\sqrt{2\pi s}}\frac{t}{t^{2}+s^{2}}ds+  \notag \\
&&+\left. \frac{1}{\pi }\frac{t}{t^{2}+s^{2}}\sum_{j=1}^{d}\frac{\partial
^{2}}{\partial x_{j}^{2}}\left( \prod_{j=1}^{d}\frac{e^{-\frac{x_{j}^{2}}{2s%
}}}{\sqrt{2\pi s}}\right) \right\vert _{0}^{+\infty },  \notag
\end{eqnarray}%
which coincides with (\ref{qua.13}).\hfill $\square $

\

\noindent \textbf{Remark 3.4 }In view of Remark 3.2, we note that the vector
process $J_{BC}^{d}(t)$ is equivalent in distribution to the following $d$%
-dimensional process:%
\begin{equation*}
\left\{
\begin{array}{c}
I_{1}(T(t)) \\
... \\
... \\
I_{d}(T(t))%
\end{array}%
\right. ,\qquad t>0,
\end{equation*}%
where $I_{j},$ $j=1,...n$ are independent iterated Brownian motions and $%
T(t) $ coincides with the first-passage time of a standard Brownian motion.

\

We focus now our attention on the process
\begin{equation*}
J_{CB}(t)=C(|B(t)|),\qquad t>0,
\end{equation*}%
and prove that its density%
\begin{equation}
p_{CB}(x,t)=\frac{2}{\pi }\int_{0}^{+\infty }\frac{s}{s^{2}+x^{2}}\frac{e^{-%
\frac{s^{2}}{2t}}}{\sqrt{2\pi t}}ds  \label{qua.15}
\end{equation}%
satisfies a non-homogeneous backward heat-equation.

\

\noindent \textbf{Theorem 3.3 }\textit{The density (\ref{qua.15}) is a
solution to the following equation}%
\begin{equation}
\frac{\partial p}{\partial t}=-\frac{1}{2}\frac{\partial ^{2}p}{\partial
x^{2}}+\frac{1}{\pi x^{2}\sqrt{2\pi t}},\qquad x\in \mathbb{R},\text{ }t>0,
\label{cr}
\end{equation}%
\textit{with initial condition }$p(x,0)=\delta (x).$

\noindent \textbf{Proof }By taking the time-derivative of $%
p_{CB}=p_{CB}(x,t) $, we have that%
\begin{eqnarray*}
\frac{\partial p}{\partial t} &=&\frac{2}{\pi }\int_{0}^{+\infty }\frac{s}{%
s^{2}+x^{2}}\frac{\partial }{\partial t}\left( \frac{e^{-\frac{s^{2}}{2t}}}{%
\sqrt{2\pi t}}\right) ds \\
&=&\frac{1}{\pi }\int_{0}^{+\infty }\frac{s}{s^{2}+x^{2}}\frac{\partial ^{2}%
}{\partial s^{2}}\left( \frac{e^{-\frac{s^{2}}{2t}}}{\sqrt{2\pi t}}\right) ds
\\
&=&\left. -\frac{1}{\pi }\int_{0}^{+\infty }\frac{\partial }{\partial s}%
\left( \frac{s}{s^{2}+x^{2}}\right) \frac{e^{-\frac{s^{2}}{2t}}}{\sqrt{2\pi t%
}}\right\vert _{0}^{+\infty }+ \\
&&+\frac{1}{\pi }\int_{0}^{+\infty }\frac{\partial ^{2}}{\partial s^{2}}%
\left( \frac{s}{s^{2}+x^{2}}\right) \frac{e^{-\frac{s^{2}}{2t}}}{\sqrt{2\pi t%
}}ds,
\end{eqnarray*}%
which coincides with (\ref{cr}).\hfill$\square $

\

The previous results can be further extended to the case where the
\textquotedblleft time process\textquotedblright\ is represented by a $n$%
-times iterated, instead of standard, Brownian motion. We prove that the
density of the process%
\begin{equation*}
C(|I_{n}(t)|)=C(|B_{1}(|B_{2}(...|B_{n+1}(t)|...)|)|),\qquad t>0
\end{equation*}%
is a solution to a non-homogeneous fractional equation.

\

\noindent \textbf{Theorem 3.4 \ }\textit{The density of the process }$%
C(|I_{n}(t)|),t>0$\textit{, which is given by}%
\begin{equation}
p_{CI}(x,t)=\frac{2}{\pi }\int_{0}^{+\infty }\frac{s}{s^{2}+x^{2}}%
p_{n}(s,t)ds,  \label{ci}
\end{equation}%
\textit{where}%
\begin{equation*}
p_{n}(x,t)=2^{n}\int_{0}^{+\infty }...\int_{0}^{+\infty }\frac{e^{-\frac{%
x^{2}}{2w_{1}}}}{\sqrt{2\pi w_{1}}}dw_{1}\frac{e^{-\frac{w_{1}^{2}}{2w_{2}}}%
}{\sqrt{2\pi w_{2}}}dw_{2}...\frac{e^{-\frac{w_{n}^{2}}{2t}}}{\sqrt{2\pi t}}%
dw_{n},
\end{equation*}%
\textit{satisfies the following equation }%
\begin{equation}
\frac{\partial ^{\frac{1}{2^{n}}}p}{\partial t^{\frac{1}{2^{n}}}}=-2^{\frac{1%
}{2^{n}}-2}\frac{\partial ^{2}p}{\partial x^{2}}+\frac{2^{n-1+\frac{1}{%
2^{n+1}}}}{\pi ^{3/2}}\frac{\Gamma \left( -\frac{1}{2}\right) }{\Gamma
\left( -\frac{1}{2^{n+1}}\right) }\frac{t^{-\frac{1}{2^{n+1}}}}{x^{2}}%
,\qquad x\in \mathbb{R},\text{ }t>0,  \label{cieq}
\end{equation}%
\textit{with initial condition }$p(x,0)=\delta (x).$

\noindent \textbf{Proof }We start by taking the fractional time-derivative
of order $1/2^{n}$ of the density (\ref{ci}):%
\begin{eqnarray}
&&\frac{\partial ^{\frac{1}{2^{n}}}}{\partial t^{\frac{1}{2^{n}}}}p_{CI}(x,t)
\label{den} \\
&=&\frac{2}{\pi }\int_{0}^{+\infty }\frac{s}{s^{2}+x^{2}}\frac{\partial ^{%
\frac{1}{2^{n}}}}{\partial t^{\frac{1}{2^{n}}}}p_{n}(s,t)ds  \notag \\
&=&\frac{2^{\frac{1}{2^{n}}-1}}{\pi }\int_{0}^{+\infty }\frac{s}{s^{2}+x^{2}}%
\frac{\partial ^{2}}{\partial s^{2}}p_{n}(s,t)ds  \notag \\
&=&\left. \frac{2^{\frac{1}{2^{n}}-1}}{\pi }\frac{s}{s^{2}+x^{2}}\frac{%
\partial }{\partial s}p_{n}(s,t)\right\vert _{0}^{+\infty }+  \notag \\
&&-\frac{2^{\frac{1}{2^{n}}-1}}{\pi }\int_{0}^{+\infty }\frac{\partial }{%
\partial s}\left( \frac{s}{s^{2}+x^{2}}\right) \frac{\partial }{\partial s}%
p_{n}(s,t)ds  \notag
\end{eqnarray}%
\begin{eqnarray*}
&=&-\left. \frac{2^{\frac{1}{2^{n}}-1}}{\pi }\frac{\partial }{\partial s}%
\left( \frac{s}{s^{2}+x^{2}}\right) p_{n}(s,t)\right\vert _{0}^{+\infty }+ \\
&&+\frac{2^{\frac{1}{2^{n}}-1}}{\pi }\int_{0}^{+\infty }\frac{\partial ^{2}}{%
\partial s^{2}}\left( \frac{s}{s^{2}+x^{2}}\right) p_{n}(s,t)ds \\
&=&\frac{2^{\frac{1}{2^{n}}-1}}{\pi }\frac{1}{x^{2}}p_{n}(0,t)-2^{\frac{1}{%
2^{n}}-2}\frac{\partial ^{2}}{\partial x^{2}}p_{CI}(x,t).
\end{eqnarray*}%
We concentrate now on the first term and evaluate%
\begin{eqnarray*}
p_{n}(0,t) &=&2^{n}\int_{0}^{+\infty }\frac{e^{-\frac{w_{1}^{2}}{2w_{2}}}}{%
\sqrt{2\pi w_{1}}}dw_{1}\int_{0}^{+\infty }\frac{e^{-\frac{w_{2}^{2}}{2w_{3}}%
}}{\sqrt{2\pi w_{2}}}dw_{2}...\int_{0}^{+\infty }\frac{e^{-\frac{w_{n}^{2}}{%
2t}}}{\sqrt{2\pi w_{n}}\sqrt{2\pi t}}dw_{n} \\
&=&\frac{2^{n-\frac{3}{4}}\Gamma \left( \frac{1}{4}\right) }{\left( \sqrt{%
2\pi }\right) ^{n+1}}\int_{0}^{+\infty }w_{2}^{-\frac{1}{4}}e^{-\frac{%
w_{2}^{2}}{2w_{3}}}dw_{2}...\int_{0}^{+\infty }\frac{e^{-\frac{w_{n}^{2}}{2t}%
}}{\sqrt{w_{n}t}}dw_{n} \\
&=&\frac{2^{n}}{\left( \sqrt{2\pi }\right) ^{n+1}}2^{-\frac{1}{2}-\frac{1}{4}%
}\Gamma \left( \frac{1}{4}\right) 2^{-\frac{1}{2}-\frac{1}{8}}\Gamma \left(
\frac{3}{8}\right) ...2^{-\frac{1}{2}-\frac{1}{2^{n}}}\Gamma \left( \frac{1}{%
2}-\frac{1}{2^{n}}\right) \int_{0}^{+\infty }\frac{w_{n}^{-\frac{1}{2^{n}}%
}e^{-\frac{w_{n}^{2}}{2t}}}{\sqrt{t}}dw_{n} \\
&=&\frac{2^{n}}{\left( \sqrt{2\pi }\right) ^{n+1}}2^{-\frac{1}{2}-\frac{1}{4}%
}\Gamma \left( \frac{1}{4}\right) 2^{-\frac{1}{2}-\frac{1}{8}}\Gamma \left(
\frac{3}{8}\right) ...2^{-\frac{1}{2}-\frac{1}{2^{n}}}\Gamma \left( \frac{1}{%
2}-\frac{1}{2^{n}}\right) 2^{-\frac{1}{2}-\frac{1}{2^{n+1}}}\Gamma \left(
\frac{1}{2}-\frac{1}{2^{n+1}}\right) t^{-\frac{1}{2^{n+1}}} \\
&=&\frac{2^{-1+\frac{1}{2^{n+1}}}}{\pi ^{\frac{n+1}{2}}}t^{-\frac{1}{2^{n+1}}%
}\Gamma \left( \frac{1}{2}-\frac{1}{2^{2}}\right) \Gamma \left( \frac{1}{2}-%
\frac{1}{2^{3}}\right) ...\Gamma \left( \frac{1}{2}-\frac{1}{2^{n}}\right)
\Gamma \left( \frac{1}{2}-\frac{1}{2^{n+1}}\right) \\
&=&\left[ \text{by the duplication property of the Gamma function}\right] \\
&=&\frac{2^{-1+\frac{1}{2^{n+1}}}}{\pi ^{\frac{n+1}{2}}}t^{-\frac{1}{2^{n+1}}%
}\pi ^{\frac{n}{2}}2^{n+1-\frac{1}{2^{n}}}\frac{\Gamma \left( -\frac{1}{2}%
\right) }{\Gamma \left( -\frac{1}{2^{n+1}}\right) } \\
&=&\frac{2^{n-\frac{1}{2^{n+1}}}}{\pi ^{\frac{1}{2}}}t^{-\frac{1}{2^{n+1}}}%
\frac{\Gamma \left( -\frac{1}{2}\right) }{\Gamma \left( -\frac{1}{2^{n+1}}%
\right) },
\end{eqnarray*}%
which, multiplied by the constant appearing in (\ref{den}), gives the final
form of equation (\ref{cieq}).\hfill$\square $

\

In the $d$-dimensional case we consider the vector process%
\begin{equation}
\left\{
\begin{array}{c}
C_{1}(|B(t)|) \\
... \\
C_{d}(|B(t)|)%
\end{array}%
\right. ,\qquad t>0  \label{vec}
\end{equation}%
and obtain the governing fractional equation in the following theorem.

\

\noindent \textbf{Theorem 3.5 \ }\textit{The joint probability law of the
process defined in (\ref{vec}) reads}%
\begin{equation}
p_{CB}^{d}(x_{1},...,x_{d},t)=\frac{2}{\pi ^{d}}\int_{0}^{+\infty
}\prod\limits_{j=1}^{d}\frac{s}{s^{2}+x_{j}^{2}}\frac{e^{-\frac{s^{2}}{2t}}}{%
\sqrt{2\pi t}}ds  \label{vec.2}
\end{equation}%
\textit{and satisfies, for }$d>1$\textit{, the following fractional equation}%
\begin{equation}
\frac{\partial p}{\partial t}=2\sum_{k=2}^{d}\sum_{j=1}^{k-1}\frac{\partial
^{2}p}{\partial |x_{k}|\partial |x_{j}|}-\sum_{k=1}^{d}\frac{\partial ^{2}p}{%
\partial x_{k}^{2}},\qquad x_{j}\in \mathbb{R},\;j=1,...,d,\text{ }t>0,
\label{tre}
\end{equation}%
\textit{with initial condition }$p(x_{1},...,x_{d},0)=\prod\limits_{j=1}^{d}%
\delta (x_{j}).$

\noindent \textbf{Proof }The time-derivative of (\ref{vec.2}) can be
evaluated as follows:%
\begin{eqnarray}
&&\frac{\partial }{\partial t}p_{CB}^{d}(x_{1},...,x_{d},t)  \label{qua} \\
&=&\frac{2}{\pi ^{d}}\int_{0}^{+\infty }\prod\limits_{j=1}^{d}\frac{s}{%
s^{2}+x_{j}^{2}}\frac{\partial }{\partial t}\left( \frac{e^{-\frac{s^{2}}{2t}%
}}{\sqrt{2\pi t}}\right) ds  \notag \\
&=&\frac{1}{\pi ^{d}}\int_{0}^{+\infty }\prod\limits_{j=1}^{d}\frac{s}{%
s^{2}+x_{j}^{2}}\frac{\partial ^{2}}{\partial s^{2}}\left( \frac{e^{-\frac{%
s^{2}}{2t}}}{\sqrt{2\pi t}}\right) ds  \notag \\
&=&\left. -\frac{1}{\pi ^{d}}\frac{\partial }{\partial s}\left(
\prod\limits_{j=1}^{d}\frac{s}{s^{2}+x_{j}^{2}}\right) \frac{e^{-\frac{s^{2}%
}{2t}}}{\sqrt{2\pi t}}\right\vert _{0}^{+\infty }+\frac{1}{\pi ^{d}}%
\int_{0}^{+\infty }\frac{\partial ^{2}}{\partial s^{2}}\left(
\prod\limits_{j=1}^{d}\frac{s}{s^{2}+x_{j}^{2}}\right) \frac{e^{-\frac{s^{2}%
}{2t}}}{\sqrt{2\pi t}}ds  \notag \\
&=&\left. -\frac{1}{\pi ^{d}}\underset{j\neq k}{\sum_{k=1}^{d}\prod%
\limits_{j=1}^{d}}\frac{s}{s^{2}+x_{j}^{2}}\frac{\partial }{\partial s}%
\left( \frac{s}{s^{2}+x_{j}^{2}}\right) \frac{e^{-\frac{s^{2}}{2t}}}{\sqrt{%
2\pi t}}\right\vert _{0}^{+\infty }+\frac{1}{\pi ^{d}}\int_{0}^{+\infty }%
\frac{\partial ^{2}}{\partial s^{2}}\left( \prod\limits_{j=1}^{d}\frac{s}{%
s^{2}+x_{j}^{2}}\right) \frac{e^{-\frac{s^{2}}{2t}}}{\sqrt{2\pi t}}ds.
\notag
\end{eqnarray}%
The first term in the last member is equal to zero (unlike the
one-dimensional case) and this makes equation (\ref{tre}) homogeneous.

Since
\begin{equation*}
\frac{\partial }{\partial s}\left( \prod\limits_{j=1}^{d}\frac{s}{%
s^{2}+x_{j}^{2}}\right) =\underset{j\neq k}{\sum_{k=1}^{d}\prod%
\limits_{j=1}^{d}}\frac{s}{s^{2}+x_{j}^{2}}\frac{\partial }{\partial s}%
\left( \frac{s}{s^{2}+x_{k}^{2}}\right) ,
\end{equation*}%
the second-order derivative becomes%
\begin{eqnarray}
\frac{\partial ^{2}}{\partial s^{2}}\left( \prod\limits_{j=1}^{d}\frac{s}{%
s^{2}+x_{j}^{2}}\right) &=&\sum_{k=1}^{d}\left[ \underset{j\neq k}{%
\sum_{j=1}^{d}}\frac{\partial }{\partial s}\left( \frac{s}{s^{2}+x_{k}^{2}}%
\right) \frac{\partial }{\partial s}\left( \frac{s}{s^{2}+x_{j}^{2}}\right)
\underset{l\neq k,j}{\prod\limits_{l=1}^{d}}\frac{s}{s^{2}+x_{l}^{2}}+\right.
\notag \\
&&\left. +\frac{\partial ^{2}}{\partial s^{2}}\left( \frac{s}{s^{2}+x_{k}^{2}%
}\right) \underset{l\neq k}{\prod\limits_{l=1}^{d}}\frac{s}{s^{2}+x_{l}^{2}}%
\right] .  \label{sei}
\end{eqnarray}%
In view of (\ref{qua.2}) and (\ref{qua.3}), we can rewrite (\ref{sei}) as%
\begin{eqnarray}
\frac{\partial ^{2}}{\partial s^{2}}\left( \prod\limits_{j=1}^{d}\frac{s}{%
s^{2}+x_{j}^{2}}\right) &=&\sum_{k=1}^{d}\left[ \underset{j\neq k}{%
\sum_{j=1}^{d}}\frac{\partial }{\partial |x_{k}|}\left( \frac{s}{%
s^{2}+x_{k}^{2}}\right) \frac{\partial }{\partial |x_{j}|}\left( \frac{s}{%
s^{2}+x_{j}^{2}}\right) \underset{l\neq k,j}{\prod\limits_{l=1}^{d}}\frac{s}{%
s^{2}+x_{l}^{2}}+\right.  \notag \\
&&\left. -\frac{\partial ^{2}}{\partial x_{k}^{2}}\left( \frac{s}{%
s^{2}+x_{k}^{2}}\right) \underset{l\neq k}{\prod\limits_{l=1}^{d}}\frac{s}{%
s^{2}+x_{l}^{2}}\right]  \label{set} \\
&=&\sum_{k=1}^{d}\underset{j\neq k}{\sum_{j=1}^{d}}\frac{\partial }{\partial
|x_{k}|}\frac{\partial }{\partial |x_{j}|}\prod\limits_{l=1}^{d}\frac{s}{%
s^{2}+x_{l}^{2}}-\sum_{k=1}^{d}\frac{\partial ^{2}}{\partial x_{k}^{2}}%
\left( \prod\limits_{l=1}^{d}\frac{s}{s^{2}+x_{l}^{2}}\right) .  \notag
\end{eqnarray}

By inserting (\ref{set}) in (\ref{qua}), we arrive at equation (\ref{tre}),
by applying the commutative property of the fractional derivative (\ref%
{qua.7}). Therefore we show that
\begin{equation*}
\frac{\partial }{\partial |x_{j}|}\frac{\partial }{\partial |x_{k}|}=\frac{%
\partial }{\partial |x_{k}|}\frac{\partial }{\partial |x_{j}|}=\frac{%
\partial ^{2}}{\partial |x_{k}|\partial |x_{j}|}.
\end{equation*}%
The second-order fractional derivative reads%
\begin{eqnarray}
&&\frac{\partial ^{2}}{\partial |y|\partial |x|}f(x,y) \\
&=&-\frac{\partial }{\partial |y|}\left\{ \frac{1}{\pi }\int_{0}^{+\infty }%
\frac{f(x-t,y)-2f(x,y)+f(x+t,y)}{t^{2}}dt\right\}  \notag \\
&=&\frac{1}{\pi ^{2}}\int_{0}^{+\infty }dz\int_{0}^{+\infty }\frac{%
f(x-t,y-z)-2f(x-t,y)+f(x-t,y+z)}{z^{2}t^{2}}dt+  \notag \\
&&-\frac{2}{\pi ^{2}}\int_{0}^{+\infty }dz\int_{0}^{+\infty }\frac{%
f(x,y-z)-2f(x,y)+f(x,y+z)}{z^{2}t^{2}}dt+  \notag \\
&&+\frac{2}{\pi ^{2}}\int_{0}^{+\infty }dz\int_{0}^{+\infty }\frac{%
f(x+t,y-z)-2f(x+t,y)+f(x+t,y+z)}{z^{2}t^{2}}dt.  \notag
\end{eqnarray}%
It is easy to check that $\frac{\partial ^{2}}{\partial |x|\partial |y|}%
f(x,y)$ produces the same result and thus the commutativity of the
second-order fractional derivative holds.\hfill $\square $

\

We consider now the last case of composition of Cauchy processes: it can be
seen that the density of the process%
\begin{equation*}
C_{1}^{a}(|C_{2}(t)|)\qquad t>0,
\end{equation*}%
where the external process is endowed with a position parameter $a\in
\mathbb{R},$ is a solution to a non-homogeneous wave equation. Indeed, by
suitably adapting the proof of theorem 4.1 in D'Ovidio and Orsingher (2009),
it is easy to check that%
\begin{equation*}
p_{CC}(x,t)=\frac{2}{\pi ^{2}}\int_{0}^{+\infty }\frac{s}{s^{2}+(x-a)^{2}}%
\frac{t}{t^{2}+s^{2}}ds
\end{equation*}%
is a solution to%
\begin{equation}
\frac{\partial ^{2}p}{\partial t^{2}}=\frac{\partial ^{2}p}{\partial x^{2}}-%
\frac{1}{\pi t(x-a)^{2}},\qquad x\in \mathbb{R},\;t>0.
\end{equation}

\

In the $d$-dimensional case the iterated Cauchy process can be defined as
follows%
\begin{equation}
J_{CC}^{d}(t)=\left\{
\begin{array}{c}
C_{1}(|C(t)|) \\
... \\
C_{d}(|C(t)|)%
\end{array}%
\right. ,\qquad t>0
\end{equation}%
where $C_{1}...C_{d}$ and $C$ are mutually independent, standard Cauchy
processes.

\

\noindent \textbf{Theorem 3.6 }\ \textit{The density of }$J_{CC}^{d}(t),t>0$%
\textit{, which can be expressed as}%
\begin{equation}
p_{CC}^{d}(x_{1},...,x_{d},t)=\frac{2}{\pi ^{d+1}}\int_{0}^{+\infty
}\prod\limits_{j=1}^{d}\frac{s}{s^{2}+x_{j}^{2}}\frac{t}{t^{2}+s^{2}}ds
\label{cc.1}
\end{equation}%
\textit{and satisfies, for }$d>1$\textit{, the following fractional equation}%
\begin{equation}
\frac{\partial ^{2}p}{\partial t^{2}}=\sum_{k=1}^{d}\frac{\partial ^{2}p}{%
\partial x_{k}^{2}}-2\sum_{k=2}^{d}\sum_{j=1}^{k-1}\frac{\partial ^{2}p}{%
\partial |x_{k}|\partial |x_{j}|},\qquad x_{j}\in \mathbb{R},\;j=1,...,d,%
\text{ }t>0,
\end{equation}%
\textit{with initial condition }$p(x_{1},...,x_{d},0)=\prod\limits_{j=1}^{d}%
\delta (x_{j}).$

\noindent \textbf{Proof \ }The second-time derivative of (\ref{cc.1}) can be
evaluated by adapting the proof of the previous theorem, as follows
\begin{eqnarray*}
&&\frac{\partial ^{2}}{\partial t^{2}}p_{CC}^{d}(x_{1},...,x_{d},t) \\
&=&\frac{2}{\pi ^{d+1}}\int_{0}^{+\infty }\prod\limits_{j=1}^{d}\frac{s}{%
s^{2}+x_{j}^{2}}\frac{\partial ^{2}}{\partial t^{2}}\left( \frac{t}{%
t^{2}+s^{2}}\right) ds \\
&=&-\frac{2}{\pi ^{d+1}}\int_{0}^{+\infty }\prod\limits_{j=1}^{d}\frac{s}{%
s^{2}+x_{j}^{2}}\frac{\partial ^{2}}{\partial s^{2}}\left( \frac{t}{%
t^{2}+s^{2}}\right) ds \\
&=&\left. \frac{2}{\pi ^{d+1}}\frac{t}{t^{2}+s^{2}}\frac{\partial }{\partial
s}\left( \prod\limits_{j=1}^{d}\frac{s}{s^{2}+x_{j}^{2}}\right) \right\vert
_{0}^{+\infty }+ \\
&&-\frac{2}{\pi ^{d+1}}\int_{0}^{+\infty }\frac{t}{t^{2}+s^{2}}\frac{%
\partial ^{2}}{\partial s^{2}}\left( \prod\limits_{j=1}^{d}\frac{s}{%
s^{2}+x_{j}^{2}}\right) ds \\
&=&\left[ \text{by (\ref{sei})}\right] \\
&=&\sum_{k=1}^{d}\frac{\partial ^{2}}{\partial x_{k}^{2}}%
p_{CC}^{d}(x_{1},...,x_{d},t)-2\sum_{k=2}^{d}\sum_{j=1}^{k-1}\frac{\partial
^{2}}{\partial |x_{k}|\partial |x_{j}|}p_{CC}^{d}(x_{1},...,x_{d},t).
\end{eqnarray*}%
\hfill$\square $

\

\noindent \textbf{Remark 3.5 }The density of $J_{CC}^{d}(t),t>0$ can be
expressed in the following alternative form%
\begin{equation}
p_{CC}^{d}(x_{1},...,x_{d},t)=\prod\limits_{j=1}^{d}\frac{2t}{\pi ^{2}\left(
t^{2}+x_{j}^{2}\right) }\ln \frac{t}{|x_{j}|},
\end{equation}%
as can be inferred from the calculations leading to theorem 4.1 of D'Ovidio
and Orsingher (2009).

\

\begin{center}
\textbf{REFERENCES}
\end{center}

\

\noindent \textbf{Allouba, H. (2002)},\textbf{\ }Brownian-time processes:
the p.d.e.connection II and the corresponding Feynman-Kac formula, \emph{%
Trans. of the Amer. Math. Soc.,} \textbf{354}, (11), 4617-4637.

\noindent \textbf{Allouba, H., Zheng, W. (2001)},\textbf{\ }Brownian-time
processes: the p.d.e.connection and half-derivative generator, \emph{%
Ann.Prob.,} \textbf{29}, (4), 1780-1795.

\noindent \textbf{Baeumer, B., Meerschaert, M., Nane, E. (2009), }Space-time
duality for fractional diffusions, \emph{arXiv: 0904.1176v1}.

\noindent \textbf{Beghin L., Orsingher E. (2003), }\textquotedblleft The
telegraph process stopped at stable-distributed times and its connection
with the fractional telegraph equation\textquotedblright , \emph{Fract.
Calc. Appl. Anal.}, \textbf{6} (2), 187-204.

\noindent \textbf{Beghin L., Orsingher E. (2009), }\textquotedblleft
Iterated elastic Brownian motions and fractional diffusion
equations\textquotedblright , \emph{Stoch. Proc. Appl.}, \textbf{119} (6);
1975-2003.

\noindent \textbf{Burdzy, K. (1994), }Variation of iterated Brownian motion,
Lecture Notes, \emph{Workshop and Conference on measure-valued processes,
stoch. partial diff. eq. and interacting syst.}, \textbf{5}, Amer. Math.
Soc., Providence, RI, 35-53.

\noindent \textbf{Chudnovsky, A., Kunin, B.} \textbf{(1987)}, A
probabilistic model of a brittle crack formation, \emph{Journ. Appl. Phys., }%
\textbf{62} (10), 4124-4129.

\noindent \textbf{De Blassie R.D. (2004), }Iterated Brownian motion in an
open set, \emph{Ann. Appl. Prob.}, \textbf{14}, (3), 1529--1558.

\noindent \textbf{D'Ovidio M., Orsingher E. (2009), }Composition of
processes and related partial differential equations. Accepted by \emph{%
Journal of Theoretical Probability}. Proofs received 16th April 2010.
Published on line 21st April 2010.

\noindent \textbf{Elmore W.C., Heald M.A. (1969)}, \emph{Physics of waves},
Dover Publ., New York.

\noindent \textbf{Khoshnevisan, D., Lewis, T.M. (1996), }A uniform modulus
result for iterated Brownian motion, \emph{Ann. Inst. Henri Poincaré},
\textbf{32} (3), 349-359.

\noindent \textbf{Khoshnevisan, D., Lewis, T.M. (1999), }Stochastic calculus
for Brownian motion on a Brownian fracture, \emph{Ann. Appl. Prob.}, \textbf{%
9} (3), 629-667.

\noindent \textbf{Nane E. (2008), }Higher-order Cauchy problems in bounded
domains, \emph{arXiv: 0809.4824v1.}

\noindent \textbf{Orsingher, E., Beghin, L. (2004),} Time-fractional
equations and telegraph processes with Brownian time, \emph{Prob. Theory and
Rel. Fields}, \textbf{128}, 141-160.

\noindent \textbf{Orsingher E., Beghin L. (2009), }Fractional diffusion
equations and processes with randomly-varying time, \emph{Ann. Prob.},
\textbf{37} (1); 206-249.

\noindent \textbf{Orsingher E., Zhao X. (1999), }Iterated processes and
their applications to higher-order differential equations, \emph{Acta Math.
Sinica}, \textbf{15} (2); 173-180.

\noindent \textbf{Podlubny, I. (1999)}, \emph{Fractional Differential
Equations}, Acad.Press, S.Diego.

\noindent \textbf{Saichev A., Zaslavsky G. (1997)},\textbf{\ }Fractional
kinetic equations: solutions and applications, \emph{Chaos,} \textbf{7},
(4), 753-764.

\

\noindent \textbf{Addresses:}

\

Luisa Beghin

Dipartimento di Statistica, Probabilità e Statistiche Applicate

\textquotedblleft Sapienza\textquotedblright\ Università di Roma

p.le A.Moro 5

00185 \ Roma \ (Italy)

e-mail: luisa.beghin@uniroma1.it

\

Lyudmyla Sakhno

Department of Mechanics and Mathematics

National Taras Shevchenko University

Kyiv, 01033 (Ukraine)

e-mail: lms@mail.univ.kiev.ua

\

Enzo Orsingher

Dipartimento di Statistica, Probabilità e Statistiche Applicate

\textquotedblleft Sapienza\textquotedblright\ Università di Roma

p.le A.Moro 5

00185 \ Roma \ (Italy)

e-mail: enzo.orsingher@uniroma1.it

\end{document}